\documentclass[11pt]{article}
\usepackage[T1]{fontenc}
\usepackage[utf8]{inputenc}
\usepackage{lmodern}
\usepackage[letterpaper,margin=1in]{geometry}
\usepackage{amsmath,amssymb,amsthm,mathtools}
\usepackage{enumitem,microtype,booktabs,xcolor}
\usepackage[colorlinks=true,linkcolor=blue!50!black,citecolor=blue!50!black,
            urlcolor=blue!50!black]{hyperref}

\hypersetup{pdftitle={Extremal adjacency and distance spectral radii of graph powers}}

\setlist[enumerate]{itemsep=3pt,topsep=4pt}
\setlist[itemize]{itemsep=3pt,topsep=4pt}
\allowdisplaybreaks[2]

\newtheorem{theorem}{Theorem}[section]
\newtheorem{lemma}[theorem]{Lemma}
\newtheorem{proposition}[theorem]{Proposition}
\newtheorem{corollary}[theorem]{Corollary}
\theoremstyle{definition}
\newtheorem{definition}[theorem]{Definition}
\newtheorem{remark}[theorem]{Remark}
\newtheorem{question}[theorem]{Question}

\newcommand{\dist}{\operatorname{dist}}
\newcommand{\diam}{\operatorname{diam}}
\newcommand{\ceil}[1]{\left\lceil #1\right\rceil}
\newcommand{\flo}[1]{\left\lfloor #1\right\rfloor}
\newcommand{\R}{\mathbb{R}}
\newcommand{\one}{\mathbf{1}}
\newcommand{\tr}{\mathsf{T}}
\newcommand{\rk}{\rho_k}
\newcommand{\muk}{\mu_k}

\title{Extremal adjacency and distance spectral radii of graph powers}
\author{Joyentanuj Das\thanks{Department of Mathematics, SRM Institute of Science and Technology Chennai, Tamil Nadu, India. Email: joyentanuj@gmail.com, joyentad@srmist.edu.in.}\ \and Anto Jerrin Ashni R\thanks{Email: rashni2001@gmail.com} \and Sangeetha C\thanks{Email: sangeethaedk12@gmail.com.}}
\date{}

\begin{document}
\maketitle

\begin{abstract}
Let $G$ be a connected graph of order $n$, and let $G^k$ denote its $k$-th power. We prove that, for every integer $k\geq1$,
\[
\rho(G^k)\geq \rho(P_n^k)
\qquad\text{and}\qquad
\rho_D(G^k)\leq \rho_D(P_n^k),
\]
where $P_n$ is the path of order $n$, and $\rho$ and $\rho_D$ are the adjacency and distance spectral radii, respectively. If $n\geq k+2$, equality in either inequality holds precisely when $G\cong P_n$. If $n\leq k+1$, every connected graph of order $n$ has complete $k$-th power, and both inequalities are equalities. The adjacency bound is proved by transferring a terminal vertex between pendant paths and applying a reflection rearrangement to a Rayleigh quotient. For the distance bound, a transfer of an entire pendant path is required; a Schur-complement comparison and estimates for the antisymmetric subspace of a path-power distance matrix yield the result. We also give examples showing why the two adjacency arguments do not extend directly to the distance spectral radius.
\end{abstract}


\section{Introduction}\label{sec:intro}

Extremal questions for graph eigenvalues have long been central to spectral graph theory; see, for example, the monograph of Brouwer and Haemers \cite{BrouwerHaemers2012}. Two quantities considered here are the adjacency spectral radius and the distance spectral radius. The former is the largest eigenvalue of the adjacency matrix, while the latter is the largest eigenvalue of the distance matrix. Among connected graphs of order $n$, the path $P_n$ is the unique graph of minimum adjacency spectral radius. This result is classical for trees, due to Collatz and Sinogowitz \cite{CollatzSinogowitz1957} and Lov\'asz and Pelik\'an \cite{LovaszPelikan1973}; the extension to connected graphs follows by taking a spanning tree and using monotonicity under edge addition. For the distance spectral radius, the extremal direction is reversed: Ruzieh and Powers \cite{RP} proved that $P_n$ is the unique maximizer among connected graphs of order $n$. Further results on distance spectra are surveyed by Aouchiche and Hansen \cite{AouchicheHansen2014}; among other results, Stevanovi\'c and Ili\'c \cite{SI} determined the unique tree of fixed order and fixed maximum degree with largest distance spectral radius.

We consider these questions for graph powers. If $k$ is a positive integer, the $k$-th power $G^k$ of a connected graph $G$ is obtained by joining every pair of distinct vertices whose distance in $G$ is at most $k$. Spectral bounds for graph powers were obtained by Das and Guo \cite{DasGuo2016}. For the second power, Chen, Li, Wang and Meng \cite{ChenLiWangMeng2019} proved a stronger $A_\alpha$-spectral result: for $0\leq\alpha<1$, the path-square has minimum $A_\alpha$-spectral radius among the squares of connected graphs of a fixed order. Thus the adjacency assertion below is already known when $k=2$. On the distance side, Xing and Zhou \cite{XZ} studied the two largest distance eigenvalues of graph powers, giving lower bounds and determining several extremal trees for distance spectra of squares. Their results do not give the path extremal inequality for every power established here. General background and bounds for the distance spectral radius can also be found in \cite{ZT2007,AouchicheHansen2014}. Transformations of pendant paths have proved useful in related extremal problems; see, for instance, \cite{NikiforovRojo2018,SI}.

For a connected graph $G$ and an integer $k\geq1$, write
\[
\rho_k(G)=\rho(G^k),\qquad \mu_k(G)=\rho_D(G^k).
\]
Our main result gives opposite extremal roles to the path for these two parameters.

\begin{theorem}\label{thm:adj}
Let $G$ be a connected graph of order $n$, and let $k\geq1$. Then
\[
\rho_k(G)\geq\rho_k(P_n).
\]
If $n\leq k+1$, equality holds for every connected graph $G$ of order $n$. If $n\geq k+2$, equality holds if and only if $G\cong P_n$.
\end{theorem}

\begin{theorem}\label{thm:dist}
Let $G$ be a connected graph of order $n$, and let $k\geq1$. Then
\[
\mu_k(G)\leq\mu_k(P_n).
\]
If $n\leq k+1$, equality holds for every connected graph $G$ of order $n$. If $n\geq k+2$, equality holds if and only if $G\cong P_n$.
\end{theorem}

The equality threshold is sharp. Indeed, if $n\leq k+1$, then every connected graph of order $n$ has diameter at most $k$, so its $k$-th power is $K_n$. For $k=1$, the two assertions reduce to the classical extremal results mentioned above. The adjacency assertion for $k=2$ is contained in the result of Chen et al. \cite{ChenLiWangMeng2019}; the present work treats all $k$ and establishes the corresponding distance inequality.

The proofs begin by reducing to trees using a spanning tree and Perron--Frobenius monotonicity. The adjacency argument uses a move that transfers one terminal vertex from a shorter pendant path to a longer one. For the distance spectral radius this move need not be monotone, and we instead transfer the whole shorter pendant path. The latter comparison is the more delicate part of the paper: after taking a Schur complement, it is reduced to a resolvent estimate for the distance matrix of a path power.

The paper is organized as follows. Section~\ref{sec:prelim} establishes notation and preliminary facts. Section~\ref{sec:adj} proves the adjacency result, including the equality characterization. Section~\ref{sec:dist} treats the distance spectral radius and gives examples explaining why the adjacency arguments cannot simply be reused. Section~\ref{sec:conc} contains concluding remarks and open problems.

\section{Preliminaries}\label{sec:prelim}

\subsection{Graphs, matrices and spectra}
All graphs are finite, simple and undirected. For a graph $H$, we write $V(H)$ and $E(H)$
for its vertex and edge sets, and $d_H(v)$ for the degree of a vertex $v$. The adjacency
matrix is denoted by $A(H)$, and its largest eigenvalue (equivalently, its spectral radius)
by $\rho(H)$. For a connected graph $H$, $\dist_H(u,v)$ denotes the distance between $u$
and $v$, and
\[
D(H)=[\dist_H(u,v)]_{u,v\in V(H)}
\]
is its distance matrix. We write $\rho_D(H)=\rho(D(H))$ for the distance spectral radius
and $\diam(H)$ for the diameter. The diagonal degree matrix of $H$ is denoted by
$\operatorname{Deg}(H)$, so that it is not confused with the distance matrix. For
$x\in V(H)$ and $r\geq0$, put
\[
B_H(x,r)=\{u\in V(H):\dist_H(u,x)\leq r\}.
\]
We write $P_n$ for the path and $K_n$ for the complete graph of order $n$. The notation
$H_1\subseteq H_2$ means that $H_1$ is a spanning subgraph of $H_2$ on the same vertex set.

Matrices are real; $\rho(M)$ is the spectral radius and, for symmetric $M$,
$\lambda_{\max}(M)=\rho(M)$ when $M\ge0$ entrywise. We use $I$ for the identity, $J$ for the
all-ones matrix, $\one$ for the all-ones vector, $e_i$ for standard basis vectors, and
$\one_S$ for the indicator vector of a set $S$; $\mathbf 1[\mathcal P]$ denotes the
indicator of a statement $\mathcal P$. For $0\le\alpha\le1$, we also use the Nikiforov matrix
\[
A_\alpha(H)=\alpha\operatorname{Deg}(H)+(1-\alpha)A(H),
\]
whose largest eigenvalue is denoted by $\lambda_\alpha(H)$; see \cite{Nikiforov2017}.
This notation is used only in the first open problem in Section~\ref{sec:conc}. Vector norms
are Euclidean and matrix norms are the induced operator norms; $\|\cdot\|_F$ is the
Frobenius norm. We use the Rayleigh--Ritz principle
\begin{equation}\label{eq:rayleigh}
\lambda_{\max}(M)=\max_{x\ne0}\frac{x^{\tr}Mx}{x^{\tr}x}\qquad(M\text{ symmetric}),
\end{equation}
see \cite[Thm.~4.2.2]{HornJohnson2013}. For a nonnegative symmetric matrix of order at least two, irreducibility is equivalent to
connectivity of its associated graph (with loops ignored); in this case its Perron vector is
strictly positive \cite[Ch.~8]{HornJohnson2013}. Distance matrices of connected graphs of
order $n\ge2$ have all off-diagonal entries positive, hence are irreducible; the one-vertex
case has $\rho=\rho_D=0$ and is always covered by the trivial range $n\le k+1$.

\subsection{Perron--Frobenius tools}
The following two lemmas are standard; we include proofs for completeness, since the exact
form of the hypotheses matters later.

\begin{lemma}\label{lem:PF}
Let $M_1,M_2$ be symmetric nonnegative matrices of the same order with $M_2$ irreducible,
$M_1\le M_2$ entrywise and $M_1\ne M_2$. Then $\rho(M_1)<\rho(M_2)$.
Consequently:
\begin{enumerate}
\item if $H_1\subseteq H_2$ are graphs on the same vertex set, $H_2$ is connected and
$H_1\ne H_2$, then $\rho(H_1)<\rho(H_2)$;
\item if $H_1,H_2$ are connected graphs on the same vertex set with
$\dist_{H_1}(u,v)\ge \dist_{H_2}(u,v)$ for all $u,v$, with strict inequality for at least one
pair, then $\rho_D(H_1)>\rho_D(H_2)$.
\end{enumerate}
\end{lemma}

\begin{proof}
Choose a nonnegative Perron vector $x$ of $M_1$ and a positive Perron vector $y$ of $M_2$.
Then
\[
(\rho(M_2)-\rho(M_1))y^{\tr}x
   =y^{\tr}(M_2-M_1)x\geq0,
\]
so $\rho(M_2)\geq\rho(M_1)$. Suppose equality holds. Since $y>0$ and
$(M_2-M_1)x\geq0$, the displayed equality forces $(M_2-M_1)x=0$. Hence
$M_2x=M_1x=\rho(M_2)x$. The irreducibility of $M_2$ now implies $x>0$. This is impossible,
since $M_2-M_1$ is a nonzero nonnegative matrix and therefore
$(M_2-M_1)x\ne0$. Thus $\rho(M_1)<\rho(M_2)$. The two graph assertions follow by applying
this matrix result to adjacency matrices and distance matrices, respectively.
\end{proof}

\begin{lemma}\label{lem:principal}
Let $M$ be a symmetric nonnegative irreducible matrix and let $P$ be a proper principal
submatrix of $M$. Then $\rho(P)<\rho(M)$.
\end{lemma}

\begin{proof}
Let $\widetilde P$ be the matrix of the same order as $M$ obtained by padding $P$ with zero
rows and columns. Then $\rho(\widetilde P)=\rho(P)$, $\widetilde P\le M$ and
$\widetilde P\ne M$ (some off-diagonal entry of $M$ in a deleted row is positive, by
irreducibility). Apply Lemma~\ref{lem:PF}.
\end{proof}

\subsection{Graph powers and distances}
For an integer $k\ge1$ define the \emph{ceiling kernel}
\[
\psi_k(t)=\ceil{\frac tk}\qquad(t\in\mathbb Z),
\]
and for a connected graph $G$ put
\[
A_k(G)=A(G^k),\qquad D_k(G)=D(G^k),\qquad \rk(G)=\rho(A_k(G)),\qquad
\muk(G)=\rho(D_k(G)).
\]
The following lemma is folklore (see, e.g., \cite{DasGuo2016,XZ}); we give the short proof.

\begin{lemma}\label{lem:distance}
For every connected graph $G$ and all $u,w\in V(G)$,
\[
\dist_{G^k}(u,w)=\psi_k\bigl(\dist_G(u,w)\bigr).
\]
\end{lemma}

\begin{proof}
A path of length $r$ in $G^k$ yields a walk of length at most $kr$ in $G$, whence
$\dist_G(u,w)\le k\,\dist_{G^k}(u,w)$ and, integrality,
$\dist_{G^k}(u,w)\ge\psi_k(\dist_G(u,w))$. Conversely, traversing a shortest $G$-path in
jumps of length at most $k$ realizes $\dist_{G^k}(u,w)\le\psi_k(\dist_G(u,w))$.
\end{proof}

\begin{corollary}\label{cor:powerfacts}
Let $G$ be connected of order $n$ and $k\ge1$. Then:
\begin{enumerate}
\item $D_k(G)=\bigl[\psi_k(\dist_G(u,w))\bigr]_{u,w}$ and $A_k(G)$ has $(u,w)$ entry
$\mathbf 1[\dist_G(u,w)\le k]$ for $u\ne w$;
\item $G^k=K_n$ if and only if $\diam(G)\le k$; in particular $n\le k+1$ implies $G^k=K_n$;
\item if $T$ is a spanning tree of $G$, then $A_k(T)\le A_k(G)$ and $D_k(T)\ge D_k(G)$
entrywise, so $\rk(T)\le\rk(G)$ and $\muk(T)\ge\muk(G)$;
\item if $T$ is a tree and $z$ a leaf of $T$ with neighbour $q$, then $T-z$ is an isometric
subtree of $T$ and $\dist_T(z,u)=1+\dist_{T-z}(q,u)$ for all $u\in V(T-z)$.
\end{enumerate}
\end{corollary}

\begin{proof}
(1) is Lemma~\ref{lem:distance}; (2) follows from (1) since $\diam(G)\le n-1$;
(3) holds because $\dist_G\le\dist_T$ pointwise and $\psi_k$ is nondecreasing, combined with
Lemma~\ref{lem:PF}; (4) is immediate since every path from $z$ starts with the edge $zq$.
\end{proof}

\subsection{Elementary properties of the ceiling kernel}

\begin{lemma}\label{lem:ceiling}
Let $k\ge1$ be an integer. Then:
\begin{enumerate}
\item $\psi_k(t+k)=\psi_k(t)+1$ for all $t\in\mathbb Z$;
\item for integers $s$ and $t\ge0$,
$\psi_k(s+t)-\psi_k(s)\ge\flo{t/k}$; if $0\le t\le k$, the difference lies in $\{0,1\}$, and
for $0<t<k$ it equals $0$ exactly when $s\bmod k\in\{1,\dots,k-t\}$;
\item $\sum_{i=0}^{k-1}\psi_k(x+i)=x+k-1$ for every $x\in\mathbb Z$.
\end{enumerate}
\end{lemma}

\begin{proof}
(1) is clear. For (2) write $t=mk+\rho$ with $0\le\rho<k$; by (1),
$\psi_k(s+t)-\psi_k(s)=m+\psi_k(s+\rho)-\psi_k(s)\ge m$, and the last difference is $0$ or
$1$, being $1$ precisely when a multiple of $k$ lies in $\{s,\dots,s+\rho-1\}$, i.e.\ when
$s\bmod k=0$ or $s\bmod k\ge k-\rho+1$. For (3) write $x=qk+\rho$, $0\le\rho<k$. If
$\rho=0$ the sum equals $kq+(k-1)$; if $\rho\ge1$ it equals
$kq+(k-\rho+1)+2(\rho-1)=kq+k+\rho-1$. In both cases it equals $x+k-1$.
\end{proof}

\subsection{Pendant paths and two grafting operations}

\begin{definition}\label{def:pendant}
A \emph{pendant path} of length $a\ge1$ at a vertex $v$ of a tree $T$ is a path
$vx_1x_2\cdots x_a$ with $\deg_T(x_i)=2$ for $1\le i\le a-1$ and $\deg_T(x_a)=1$; we set
$x_0=v$. A vertex of degree at least $3$ is a \emph{branching vertex}. A tree with at most one
branching vertex is a \emph{spider}; $S(\ell_1,\dots,\ell_t)$ denotes the spider with $t\ge2$
legs (pendant paths) of lengths $\ell_1,\dots,\ell_t\ge1$ at its centre, a tree of order
$1+\sum_i\ell_i$. A tree with exactly three leaves is a spider $S(\ell_1,\ell_2,\ell_3)$.
\end{definition}

\begin{definition}\label{def:moves}
Let $T$ be a tree with two pendant paths $vx_1\cdots x_a$ and $vy_1\cdots y_b$ at a common
vertex $v$, where $a\ge b\ge1$, and put $y_0=v$.
\begin{itemize}
\item[\rm(M1)] \emph{Terminal-leaf move.} $T'$ is obtained from $T$ by deleting the edge
$y_{b-1}y_b$ and adding the edge $x_ay_b$.
\item[\rm(M2)] \emph{Whole-path graft.} $T^\star$ is obtained from $T$ by deleting the edge
$vy_1$ and adding the edge $x_ay_1$.
\end{itemize}
Both $T'$ and $T^\star$ are trees of the same order as $T$. In (M2) we further let $C$ be the
component containing $v$ of $T-\{x_1,\dots,x_a,y_1,\dots,y_b\}$, let
$U=V(C)\setminus\{v\}$, $h_u=\dist_C(u,v)$ and
\[
H=\max_{u\in V(C)}\dist_C(u,v).
\]
Thus $C=\{v\}$ exactly when $T$ is a path, and $H\ge1$ otherwise.
\end{definition}

\begin{remark}\label{rem:composition}
Up to isomorphism, (M2) is the composition of $b$ moves of type (M1). Indeed, after $i$ moves
($0\le i<b$) the two pendant paths at $v$ have lengths $a+i$ and $b-i$, and $a+i\ge b-i$; the
next (M1) move transfers the terminal vertex of the shorter one to the end of the longer one.
After $b$ steps the shorter path has disappeared and the longer one has length $a+b$, which
is $T^\star$ up to reversing the labels of $Y$.
\end{remark}

\section{The adjacency spectral radius of graph powers}\label{sec:adj}

\subsection{A hanging-path lemma for powers}

\begin{lemma}\label{lem:hanging}
Let $T'$ be obtained from $T$ by the move \rm(M1) of Definition~\ref{def:moves}. Then
\[
\rk(T')\ \le\ \rk(T)\qquad\text{for every integer }k\ge1 .
\]
\end{lemma}

\begin{proof}
Put $z=y_b$, $T^-=T-z$, and
\[
q_0=y_{b-1},\qquad q_L=x_a,\qquad L=a+b-1 .
\]
The unique $q_0$--$q_L$ path in $T^-$ is $Q=q_0q_1\cdots q_L$, where $q_i=y_{b-1-i}$ for
$0\le i\le b-1$ and $q_{b-1+i}=x_i$ for $1\le i\le a$. Writing $s=b-1$ we have
$L-2s=a-b+1\ge1$, so $s<L/2$; moreover every component of $T^--V(Q)$ is attached at
$v=q_s$, because all internal vertices of the two pendant paths have degree $2$ in $T$.

By Corollary~\ref{cor:powerfacts}(4), in $T$ the vertex $z$ is adjacent in $T^k$ exactly to
$S_0:=B_{T^-}(q_0,r)$, and in $T'$ exactly to $S_1:=B_{T^-}(q_L,r)$, where $r=k-1$. Also, the
subgraph of either power induced on $V(T^-)$ is $(T^-)^k$. Hence, with $A=A_k(T^-)$ and in the
vertex order $(z,V(T^-))$,
\[
A_k(T)=\begin{pmatrix}0&\one_{S_0}^{\tr}\\ \one_{S_0}&A\end{pmatrix},
\qquad
A_k(T')=\begin{pmatrix}0&\one_{S_1}^{\tr}\\ \one_{S_1}&A\end{pmatrix}.
\]

\paragraph{The compression.}
For $\xi\in\R_{\ge0}^{V(T^-)}$ define $\mathcal C\xi$ by sorting each mirror pair
$\{q_i,q_{L-i}\}$, $0\le i<L/2$, so that the larger coordinate sits on the left:
\[
(\mathcal C\xi)(q_i)=\max\{\xi(q_i),\xi(q_{L-i})\},\qquad
(\mathcal C\xi)(q_{L-i})=\min\{\xi(q_i),\xi(q_{L-i})\}.
\]
The midpoint $q_{L/2}$ (if $L$ is even) and all vertices outside $Q$ are fixed. Clearly
$\|\mathcal C\xi\|=\|\xi\|$. Write $\alpha_i=\xi(q_i)$, $\beta_i=\xi(q_{L-i})$,
$d_i=\alpha_i-\beta_i$, and let $d_i'=|d_i|$ be the corresponding quantity for
$\mathcal C\xi$; the pair sums $\alpha_i+\beta_i$ are preserved.

\paragraph{Claim 1 (energy).} $(\mathcal C\xi)^{\tr}A(\mathcal C\xi)\ge\xi^{\tr}A\xi$.

Edges inside one mirror pair, edges between two vertices outside $Q$, and edges between a
fixed vertex and a mirror pair member other than through the pair structure are unchanged or
handled below; we consider the three types of contributions.

\emph{(a) Two mirror pairs.} Let $0\le i<j\le L/2$ (for $j=L/2$ the ``pair'' is the fixed
midpoint). The parallel distances are $\dist(q_i,q_j)=\dist(q_{L-i},q_{L-j})=j-i$, while the
cross distances are $\dist(q_i,q_{L-j})=\dist(q_j,q_{L-i})=L-i-j$. Since $2j\le L$ we have
$j-i\le L-i-j$; hence a cross edge forces both parallel edges. Writing
$p=\mathbf 1[j-i\le k]$ and $c=\mathbf 1[L-i-j\le k]$ (so $c\le p$), the contribution of these
two pairs to the quadratic form equals
\[
2p(\alpha_i\alpha_j+\beta_i\beta_j)+2c(\alpha_i\beta_j+\alpha_j\beta_i)
=(p+c)(\alpha_i+\beta_i)(\alpha_j+\beta_j)+(p-c)\,d_id_j .
\]
The first term is invariant under $\mathcal C$, and the second changes by
$(p-c)(|d_id_j|-d_id_j)\ge0$. (For $j=L/2$ one has $p=c$, so nothing changes.)

\emph{(b) A vertex outside $Q$ and a mirror pair.} Let $w\notin Q$ lie in a component attached
at $q_s$, and put $h=\dist(w,q_s)$. For $i\le L/2$,
\[
\dist(w,q_i)=h+|i-s|\ \le\ h+L-i-s=\dist(w,q_{L-i}),
\]
because $|i-s|\le L-i-s$ follows from $2i\le L$ when $i\ge s$ and from $2s\le L$ when $i<s$.
Hence $w$ is adjacent to $q_{L-i}$ only if it is adjacent to $q_i$. If both adjacencies are
present or both absent, the contribution is $\xi(w)(\alpha_i+\beta_i)$ or $0$, unchanged; if
only the left one is present, the contribution rises from $\xi(w)\alpha_i$ to
$\xi(w)\max\{\alpha_i,\beta_i\}$.

\emph{(c) Remaining edges.} Edges within a mirror pair contribute $2\alpha_i\beta_i$, which is
invariant; edges with both ends outside $Q$, or joining an outside vertex to the midpoint,
are untouched. This proves Claim~1.

\paragraph{Claim 2 (balls).} $\one_{S_0}^{\tr}\mathcal C\xi\ \ge\ \one_{S_1}^{\tr}\xi$.

For a mirror pair with $i<L/2$ we have $\dist(q_i,q_0)=\dist(q_{L-i},q_L)=i$ and
$\dist(q_{L-i},q_0)=\dist(q_i,q_L)=L-i$. If $i\le r<L-i$, the left side receives
$\max\{\alpha_i,\beta_i\}$ and the right side $\beta_i$. If $L-i\le r$, both balls contain the
whole pair and the pair sum is preserved. If $r<i$, neither ball meets the pair. The midpoint
contributes equally to both sides. Finally, for $w\notin Q$ attached at $q_s$ with
$s\le L/2$, $\dist(w,q_0)=h+s\le h+L-s=\dist(w,q_L)$, so $w\in S_1$ implies $w\in S_0$; since
$\mathcal C$ fixes $w$, the off-path contribution to the left is at least that to the right.
This proves Claim~2.

\paragraph{Conclusion.} Let $f=(\alpha,\xi)>0$ be a Perron vector of the connected graph
$(T')^k$, with $\alpha=f(z)$. Then
\[
\rk(T')\|f\|^2=\xi^{\tr}A\xi+2\alpha\,\one_{S_1}^{\tr}\xi
\ \le\ (\mathcal C\xi)^{\tr}A(\mathcal C\xi)+2\alpha\,\one_{S_0}^{\tr}\mathcal C\xi
=(\alpha,\mathcal C\xi)^{\tr}A_k(T)(\alpha,\mathcal C\xi),
\]
using Claims 1 and 2. Since $\|(\alpha,\mathcal C\xi)\|=\|f\|$, \eqref{eq:rayleigh} gives
$\rk(T')\|f\|^2\le\rk(T)\|f\|^2$.
\end{proof}

\subsection{Strictness of the terminal move}

We first record a resolvent comparison for paths, used for $k=1$. Let $A_m=A(P_m)$ and
\[
\Delta_m(t)=\det(tI-A_m)\quad(m\ge0),\qquad \Delta_{-1}(t)=0,\ \Delta_0(t)=1 ,
\]
so that $\Delta_m=t\Delta_{m-1}-\Delta_{m-2}$ for $m\ge1$. The eigenvalues of $P_m$ are
$2\cos\frac{\pi j}{m+1}$, $1\le j\le m$ \cite[\S1.4.3]{BrouwerHaemers2012}; hence
$\Delta_m(t)>0$ for all $0\le m\le N$ whenever $t>\rho(P_N)$.

\begin{lemma}\label{lem:resolvent-path}
Let $N\ge3$ and $2\le j\le N-1$. For every $t>\rho(P_N)$,
\[
\bigl((tI-A_N)^{-1}\bigr)_{j,j}\ >\ \bigl((tI-A_N)^{-1}\bigr)_{N,N}.
\]
\end{lemma}

\begin{proof}
Deleting vertex $j$ from $P_N$ leaves $P_{j-1}\cup P_{N-j}$, so by the cofactor formula
\[
\bigl((tI-A_N)^{-1}\bigr)_{j,j}=\frac{\Delta_{j-1}(t)\Delta_{N-j}(t)}{\Delta_N(t)},
\qquad
\bigl((tI-A_N)^{-1}\bigr)_{N,N}=\frac{\Delta_{N-1}(t)}{\Delta_N(t)} .
\]
Put $D(a,b)=\Delta_a\Delta_b-\Delta_{a-1}\Delta_{b-1}$ for $a,b\ge0$. The recurrence in $b$
gives $D(a,b)=tD(a,b-1)-D(a,b-2)$ for $b\ge2$, with $D(a,0)=\Delta_a$ and
$D(a,1)=t\Delta_a-\Delta_{a-1}=\Delta_{a+1}$. Since $\Delta_{a+b}$ obeys the same recurrence
and initial values in $b$, we get the identity
\[
\Delta_a(t)\Delta_b(t)-\Delta_{a-1}(t)\Delta_{b-1}(t)=\Delta_{a+b}(t)\qquad(a,b\ge0).
\]
With $a=j-1\ge1$ and $b=N-j\ge1$ this yields
$\Delta_{j-1}\Delta_{N-j}=\Delta_{N-1}+\Delta_{j-2}\Delta_{N-j-1}>\Delta_{N-1}$ for
$t>\rho(P_N)$, since all factors are positive. Divide by $\Delta_N(t)>0$.
\end{proof}

\begin{lemma}\label{lem:terminal-strict}
Let $k\ge1$ and $n\ge k+2$, and let $U$ be a non-path tree of order $n$ that is carried to
$P_n$ by a single move \rm(M1). Then
\[
\rk(U)\ >\ \rk(P_n).
\]
\end{lemma}

\begin{proof}
\emph{Structure.} Put $z=y_b$ and $T^-=U-z$. Since the resulting tree is $P_n$ and $z$
is an end-vertex of that path, $T^-$ is a path on $n-1$ vertices, with $x_a$ as one
endpoint. If $b\ge2$, then $y_{b-1}$ has degree $1$ in $T^-$ and is its other endpoint.
Attaching $z$ to this endpoint would leave a path, contrary to the assumption that $U$ is
not a path. Thus $b=1$ and $y_{b-1}=v$. The vertex $v$ is internal in $T^-$; otherwise
attaching $z$ to $v$ would again give a path. Moreover, $v$ has degree $3$ in $U$, so $U$
is a spider with three legs. One leg of $T^-$ runs from $v$ to $q_L=x_a$ and has length
$L=\dist_{T^-}(v,q_L)=a$; the other runs from $v$ to the opposite endpoint and has some
length $c\ge1$. Consequently $L+c=n-2\geq k$.

\emph{Case $k\ge2$.} Put $q_0=v$, $r=k-1$, $S_0=B_{T^-}(q_0,r)$, $S_1=B_{T^-}(q_L,r)$ and let $Q$ be the
$q_0$--$q_L$ path. All vertices of $T^--V(Q)$ hang at the left endpoint $q_0$ of $Q$, i.e.\ at
$q_s$ with $s=0\le L/2$; hence the compression $\mathcal C$ of Lemma~\ref{lem:hanging}
satisfies Claims 1 and 2 there:
\[
(\mathcal C\xi)^{\tr}A_k(T^-)(\mathcal C\xi)\ge\xi^{\tr}A_k(T^-)\xi,
\qquad \one_{S_0}^{\tr}\mathcal C\xi\ge\one_{S_1}^{\tr}\xi .
\]
Choose $w$ on the side of $T^-$ opposite to $q_L$ at distance $d=\max\{1,k-L\}$ from $q_0$.
This vertex exists: if $k-L\ge1$ then $c=n-2-L\ge k-L=d$, and otherwise $d=1\le c$. Also
$d\le k-1=r$, because $L\ge1$ and $k\ge2$. Hence $w\in S_0$, whereas
$\dist_{T^-}(w,q_L)=d+L\ge k>r$, so $w\notin S_1$.

Let $f=(\alpha,\xi)>0$ be the Perron vector of $P_n^k$, with $\alpha=f(z)$; note that $z$ is
attached to $q_L$ in $P_n$. Since $\mathcal C$ fixes $w\notin Q$ and $\xi(w)>0$, Claim 2 is
\emph{strict}: $\one_{S_0}^{\tr}\mathcal C\xi>\one_{S_1}^{\tr}\xi$. Therefore, by
\eqref{eq:rayleigh},
\[
\rk(U)\ \ge\ \frac{(\alpha,\mathcal C\xi)^{\tr}A_k(U)(\alpha,\mathcal C\xi)}{\|(\alpha,\mathcal C\xi)\|^2}
\ >\ \frac{(\alpha,\xi)^{\tr}A_k(P_n)(\alpha,\xi)}{\|(\alpha,\xi)\|^2}\ =\ \rk(P_n).
\]

\emph{Case $k=1$.} A non-path tree of this form has $n\ge4$. Put $R=P_N$, where
$N=n-1\ge3$. Then $U$ is obtained by attaching a leaf $z$ to
an internal vertex $j$ ($2\le j\le N-1$), while $P_n$ arises by attaching $z$ to the endpoint
$N$. For $s\in\{j,N\}$ let $H_s$ be the corresponding graph and $\lambda_s=\rho(H_s)$. By
Lemma~\ref{lem:PF}, $\lambda_s>\rho(P_N)$, and the Perron equations for
$A(H_s)=\begin{pmatrix}0&e_s^{\tr}\\ e_s&A_N\end{pmatrix}$ give
$\lambda_s=e_s^{\tr}(\lambda_sI-A_N)^{-1}e_s$. The function
$\phi_s(t)=e_s^{\tr}(tI-A_N)^{-1}e_s-t$ is strictly decreasing on $(\rho(P_N),\infty)$, since
$\phi_s'(t)=-\|(tI-A_N)^{-1}e_s\|^2-1<0$, and $\phi_s(\lambda_s)=0$. By
Lemma~\ref{lem:resolvent-path}, $\phi_j(t)>\phi_N(t)$ there; in particular
$\phi_j(\lambda_N)>\phi_N(\lambda_N)=0$, whence $\lambda_j>\lambda_N$, i.e.\
$\rho(U)>\rho(P_n)$.
\end{proof}

\subsection{Trees, and the proof of Theorem~\ref{thm:adj}}

\begin{proposition}\label{prop:trees-adj}
Let $T$ be a tree of order $n$ and $k\ge1$. Then $\rk(T)\ge\rk(P_n)$. If $n\ge k+2$, equality
holds if and only if $T\cong P_n$.
\end{proposition}

\begin{proof}
If $T\cong P_n$ there is nothing to prove. Otherwise root $T$ at a leaf and let $v$ be a
branching vertex farthest from the root. Every component of $T-v$ not containing the root is
a pendant path, and there are at least two of them, say of lengths $a\ge b\ge1$ with $a$
maximal. By Lemma~\ref{lem:hanging}, the move (M1) applied to these two paths produces a tree
$T'$ with $\rk(T')\le\rk(T)$. No new branching vertex is created ($x_a$ acquires degree $2$,
$y_{b-1}$ loses one neighbour), and the total length of the downward pendant paths at $v$
other than a longest one drops by exactly $1$. Hence after finitely many such moves $v$ has a
single downward pendant path and degree $2$, so the number of branching vertices has
decreased. Iterating, we obtain a finite sequence $T=T_0,T_1,\dots,T_m=P_n$ with
$\rk(T_{i+1})\le\rk(T_i)$, which proves the inequality.

Assume $n\ge k+2$ and $T\not\cong P_n$. The sequence has at least one step, and its last step
transforms a non-path tree into $P_n$; by Lemma~\ref{lem:terminal-strict} that step is
strict. Hence $\rk(T)>\rk(P_n)$.
\end{proof}

\begin{proof}[Proof of Theorem~\ref{thm:adj}]
Let $T$ be a spanning tree of $G$. By Corollary~\ref{cor:powerfacts}(3) and
Proposition~\ref{prop:trees-adj},
\[
\rk(G)\ \ge\ \rk(T)\ \ge\ \rk(P_n).
\]
If $n\le k+1$, then $G^k=K_n=P_n^k$ by Corollary~\ref{cor:powerfacts}(2), so equality holds
and is equivalent to $G^k\cong P_n^k$.

Assume $n\ge k+2$ and $\rk(G)=\rk(P_n)$. Then $\rk(G)=\rk(T)=\rk(P_n)$, so
Proposition~\ref{prop:trees-adj} forces $T\cong P_n$, and Lemma~\ref{lem:PF}(1) applied to
$T^k\subseteq G^k$ forces $T^k=G^k$. Label the spanning path $T$ as $1,2,\dots,n$, so
$G^k=P_n^k$.

Suppose $G\ne P_n$; then $G$ contains an edge $ij$ with $s:=j-i\ge2$. If $s\ge k+1$ then
$ij\in E(G^k)\setminus E(P_n^k)$, a contradiction. So $2\le s\le k$. Since $n\ge k+2$,
\[
(i-1)+(n-j)=n-s-1\ \ge\ k+1-s\ \ge\ 1 ,
\]
so we may choose integers $t,u\ge0$ with $t\le i-1$, $u\le n-j$ and $t+u=k+1-s$. Put
$p=i-t\ge1$ and $q=j+u\le n$. Along the path, $q-p=t+s+u=k+1$, so $pq\notin E(P_n^k)$. But in
$G$ we may travel from $p$ to $i$, use the edge $ij$, and travel from $j$ to $q$, whence
\[
\dist_G(p,q)\le t+1+u=k+2-s\le k ,
\]
so $pq\in E(G^k)$, contradicting $G^k=P_n^k$. Therefore $G\cong P_n$. Conversely, $G^k\cong
P_n^k$ trivially implies equality, and $G\cong P_n$ implies $G^k\cong P_n^k$.
\end{proof}

\begin{corollary}\label{cor:graft-adj}
If $T^\star$ is obtained from $T$ by the move \rm(M2), then $\rk(T^\star)\le\rk(T)$ for every
$k\ge1$.
\end{corollary}

\begin{proof}
By Remark~\ref{rem:composition}, $T^\star$ is obtained from $T$ by $b$ moves of type (M1),
each satisfying the hypothesis $a\ge b$ of Lemma~\ref{lem:hanging}, which we apply $b$ times.
\end{proof}

\begin{corollary}\label{cor:sharp-adj}
If $n\le k+1$, every connected graph $G$ of order $n$ satisfies $G^k=K_n=P_n^k$, so all
connected graphs of order $n$ are extremal in Theorem~\ref{thm:adj}. For $k\ge2$ the condition
$n\ge k+2$ is therefore sharp.
\end{corollary}

\begin{proof}
Immediate from Corollary~\ref{cor:powerfacts}(2); for $k\ge2$ and $n=k+1\ge3$ there exist
connected $n$-vertex graphs other than $P_n$.
\end{proof}

\section{The distance spectral radius of graph powers}\label{sec:dist}

\subsection{Spectral estimates for path-power distance matrices}

Throughout this subsection fix integers $a\ge1$, $k\ge1$, put $N=a+1$, and index the vertices
of a path by $0,1,\dots,a$. Let
\[
P=P_{a,k}=\bigl[\psi_k(|i-j|)\bigr]_{i,j=0}^{a}=D_k(P_{a+1}),
\qquad
\bar\sigma=\frac{\one^{\tr}P\one}{N},
\]
so $\bar\sigma$ is the average row sum of $P$ and $\rho(P)\ge\bar\sigma$. Let $R$ be the
reversal permutation matrix, $(Rx)_i=x_{a-i}$, and let
\[
\mathcal A=\{x\in\R^N:\ Rx=-x\}
\]
be the antisymmetric subspace, of dimension $\flo{N/2}$. Since $P$ commutes with $R$, both
$\mathcal A$ and its orthogonal complement are $P$-invariant; let
$\gamma=\lambda_{\max}(P|_{\mathcal A})$. For integers $t\ge1$ put
\begin{equation}\label{eq:qdef}
c_t=\bigl(\psi_k(i+t)\bigr)_{i=0}^{a},\qquad
q_t=(R-I)c_t,\qquad q_t(i)=\psi_k(a-i+t)-\psi_k(i+t).
\end{equation}
Then $q_t\in\mathcal A$, and $q_{t+k}=q_t$ by Lemma~\ref{lem:ceiling}(1). Finally, for
$\lambda>\rho(P)$ let $K=(\lambda I-P)^{-1}$; $K$ is positive definite, commutes with $R$,
and $K|_{\mathcal A}$ is positive definite with operator norm $(\lambda-\gamma)^{-1}$.
Consequently
\begin{equation}\label{eq:resolvent}
\Bigl|\tfrac12 q_j^{\tr}Kq_h\Bigr|\ \le\ \frac{\|q_j\|\,\|q_h\|}{2(\lambda-\gamma)}
\qquad(j,h\ge1).
\end{equation}
The aim of this subsection is Corollary~\ref{cor:resolvent}: the left-hand side of
\eqref{eq:resolvent} is \emph{strictly smaller} than $\flo{a/k}$, and hence smaller than the
ceiling increment $\psi_k(a+j+h)-\psi_k(j+h)$ (Lemma~\ref{lem:ceiling}(2)). We treat the
ranges $k\le a<2k$ and $a\ge2k$ separately, and use a uniform norm bound in both.

\begin{lemma}\label{lem:middle}
Let $a=k+r$ with $0\le r\le k-1$. Then
\begin{align}
\bar\sigma&=a+\frac{r(r+1)}{N},\label{eq:middleL}\\
\gamma&\le\frac r2-1,\label{eq:middlegamma}\\
\|q_t\|^2&\le\max\{N,\,4(r+1)\}\qquad(t\ge1).\label{eq:middlenorm}
\end{align}
Moreover, with $g=N+\dfrac{r(r+1)}N-\dfrac r2$ we have $\lambda-\gamma>g$ for every
$\lambda>\rho(P)$, and
\begin{equation}\label{eq:middlecompare}
\max\{N,4(r+1)\}\ \le\ 2g .
\end{equation}
\end{lemma}

\begin{proof}
Here $\psi_k(d)\in\{1,2\}$ for $1\le d\le a$, so $P=J-I+B$ with
$B_{ij}=\mathbf 1[|i-j|>k]$. The number of unordered pairs at distance greater than $k$ is
$\sum_{d=k+1}^{k+r}(N-d)=\sum_{j=1}^{r}j=\frac{r(r+1)}2$, giving
$\one^{\tr}P\one=N(N-1)+r(r+1)$ and hence \eqref{eq:middleL}.

In the orthonormal antisymmetric basis $u_i=(e_i-e_{a-i})/\sqrt2$, $0\le i<\flo{N/2}$, we
have $u_i^{\tr}Bu_j=B_{ij}-B_{i,a-j}$. Since $i,j<N/2\le k$, $B_{ij}=0$; and
$B_{i,a-j}=1$ iff $a-i-j>k$, i.e.\ $i+j\le r-1$. Hence $B|_{\mathcal A}=-F_r\oplus0$, where
$F_r=\bigl(\mathbf 1[i+j\le r-1]\bigr)_{i,j=0}^{r-1}$ (and $F_0$ is empty); note
$N=k+r+1\ge2r+2$, so $\flo{N/2}\ge r+1$ and $F_r$ indeed fits. Write
$F_r=\frac12J_r+E_r$; every entry of $E_r$ is $\pm\frac12$, so
$\|E_r\|\le\|E_r\|_F=\frac r2$ and $\lambda_{\min}(F_r)\ge-\frac r2$. Thus
$\lambda_{\max}(B|_{\mathcal A})\le\max\{\frac r2,0\}=\frac r2$. As $Jx=0$ for
$x\in\mathcal A$, $(J-I)|_{\mathcal A}=-I$, which proves \eqref{eq:middlegamma}.

For \eqref{eq:middlenorm} we may assume $1\le t\le k$ by periodicity. Let $0\le i<a/2$; then
$i\le k-1$. Set $A=r+t-k$ and $B_0=k-t+1$, so $A+B_0=r+1$. Since $1\le i+t\le2k-1$,
$\psi_k(i+t)=1+\mathbf 1[i\ge B_0]$; and since $a-i+t=2k+A-i$ with $A-i\in(-2k,k)$,
$\psi_k(a-i+t)=2+\mathbf 1[i<A]-\mathbf 1[i\ge A+k]$. Subtracting,
\begin{equation}\label{eq:cutoffs}
q_t(i)=1+\mathbf 1[i<A]-\mathbf 1[i\ge A+k]-\mathbf 1[i\ge B_0].
\end{equation}
If $A<0$, then $q_t(i)\le1$; also $q_t(i)\ge0$ on the left half, since $a-i+t\ge i+t$ there.
Hence $q_t(i)\in\{0,1\}$, and antisymmetry gives $|q_t(i)|\le1$ for all $i$, so
$\|q_t\|^2\le N$.

If $A\ge0$, then $i<k\le A+k$ and \eqref{eq:cutoffs} becomes
$q_t(i)=\mathbf 1[i<A]+\mathbf 1[i<B_0]$ on the left half. Moreover $A\ge0$ forces
$t\ge k-r$, hence $B_0\le r+1$, while $A\le r$; both cut-offs are at most
$r+1\le\flo{N/2}$, i.e.\ they lie inside the left half. Summing squares over the left half
and reflecting (a central coordinate, if present, is $0$),
\[
\|q_t\|^2=2\bigl(A+B_0+2\min(A,B_0)\bigr)\le4(A+B_0)=4(r+1).
\]
Finally, $\lambda>\rho(P)\ge\bar\sigma$ together with \eqref{eq:middleL} and
\eqref{eq:middlegamma} gives $\lambda-\gamma>g$. Clearly $2g\ge 2N-r\ge N$. For the other
bound put $d=k-r-1\ge0$, so $N=2r+2+d$; a direct simplification gives
\[
2g-4(r+1)=2d-r+\frac{2r(r+1)}N=d\Bigl(2-\frac rN\Bigr)\ \ge\ 0 ,
\]
which proves \eqref{eq:middlecompare}.
\end{proof}

\begin{lemma}\label{lem:norm}
For all integers $a\ge0$, $k\ge1$ and $t\ge1$,
\begin{equation}\label{eq:norm}
\|q_t\|^2\ \le\ (a+1)\Bigl(1+\frac{a(a+2)}{3k^2}\Bigr).
\end{equation}
\end{lemma}

\begin{proof}
Let $\ell_a(i)=(a-2i)/k$; note $\ell_a(a-i)=-\ell_a(i)$ and
$\sum_{i=0}^{a}\ell_a(i)^2=\frac{a(a+1)(a+2)}{3k^2}$. We prove
\begin{equation}\label{eq:excess}
E(a):=\sum_{i=0}^{a}\bigl(q_t(i)^2-\ell_a(i)^2\bigr)\ \le\ a+1
\end{equation}
by induction on $a$; \eqref{eq:norm} then follows.

\emph{Base $a<k$.} The numbers $\psi_k(i+t)$, $0\le i\le a$, lie in an interval of length
$a<k$, hence take at most two consecutive integer values, so $|q_t(i)|\le1$ and
$E(a)\le a+1$.

\emph{Base $k\le a<2k$.} Write $a=k+r$; by Lemma~\ref{lem:middle} it suffices to check
$4(r+1)\le(a+1)\bigl(1+\frac{a(a+2)}{3k^2}\bigr)$, the bound $N$ being trivially admissible.
Put $s=r+1\ge1$ and $d=k-r-1\ge0$, so $k=s+d$ and $a+1=2s+d$. Multiplying by $3k^2$, the
difference between right and left sides equals
\[
2s(s^2-1)+d(3s^2-1)+6sd^2+4d^3\ \ge\ 0 .
\]

\emph{Step $a\ge2k$.} Delete the first and last $k$ coordinates. For $i=i'+k$ with
$0\le i'\le a-2k$, Lemma~\ref{lem:ceiling}(1) gives
$q_t(i)=\psi_k\bigl((a-2k)-i'+t\bigr)-\psi_k(i'+t)$ and $\ell_a(i)=\ell_{a-2k}(i')$; i.e.\
the middle block is exactly the pair of vectors with parameter $a-2k$, contributing at most
$(a-2k)+1$ by induction.

On the first block $0\le i\le k-1$, the values $q_t(i)$ lie in an interval of length at most
two: as $i$ runs over $k$ consecutive integers, $\psi_k(a-i+t)$ drops by at most one and
$\psi_k(i+t)$ rises by at most one. By Lemma~\ref{lem:ceiling}(3),
\[
\sum_{i=0}^{k-1}q_t(i)=(a+t)-(t+k-1)=a-k+1=\sum_{i=0}^{k-1}\ell_a(i),
\]
so both blocks have the same mean $\mu$. Numbers in an interval of length two have variance
at most one (Popoviciu's bound: variance $\le(\text{range})^2/4$), while variance is
nonnegative; hence
\[
\sum_{i=0}^{k-1}q_t(i)^2-\sum_{i=0}^{k-1}\ell_a(i)^2
\ \le\ \bigl(k+k\mu^2\bigr)-k\mu^2=k .
\]
The last block obeys the same bound by antisymmetry of $q_t$ and $\ell_a$. Adding the three
contributions gives $E(a)\le(a-2k+1)+2k=a+1$.
\end{proof}

\begin{lemma}\label{lem:long}
Let $k\ge2$ and $a\ge2k$, and put
\[
\eta=\frac{k-1}{2k},\qquad g=\frac{a(a+2)}{3k}+\eta,\qquad m=\flo{\frac ak}\ \ge2 .
\]
Then $\lambda-\gamma>g$ for every $\lambda>\rho(P)$, and
\begin{equation}\label{eq:longnormgap}
(a+1)\Bigl(1+\frac{a(a+2)}{3k^2}\Bigr)\ <\ 2mg .
\end{equation}
\end{lemma}

\begin{proof}
Let $D=[\,|i-j|\,]_{i,j=0}^{a}$. If $x^{\tr}\one=0$, then splitting along the $a$ cuts of the
path and writing $S_r=\sum_{i\le r}x_i$ gives
\begin{equation}\label{eq:pathnegative}
x^{\tr}Dx=-2\sum_{r=0}^{a-1}S_r^2\ \le\ 0 .
\end{equation}
Write $P=\frac1kD+\eta(J-I)+E$. Indeed, for $i\ne j$ we have
$P_{ij}=\frac{|i-j|}k+\varepsilon_{|i-j|}$ with $\varepsilon_d=\psi_k(d)-d/k\in[0,2\eta]$,
so $E$ has zero diagonal and $|E_{ij}|\le\eta$; hence
$\|E\|\le\max_i\sum_j|E_{ij}|\le a\eta$. For $x\in\mathcal A$ we have $x^{\tr}\one=0$, so
\eqref{eq:pathnegative} yields
$x^{\tr}Px\le-\eta\|x\|^2+a\eta\|x\|^2$, i.e.
\begin{equation}\label{eq:gammalong}
\gamma\ \le\ (a-1)\eta .
\end{equation}
Next, $\varepsilon_d$ has period $k$ and mean $\eta$, and for $p\equiv\ell\pmod k$,
$0\le\ell<k$, a direct computation gives
$\sum_{d=1}^{p}(\varepsilon_d-\eta)=\frac{\ell(k-\ell)}{2k}\ge0$. Since the weights
$N-d$ are positive and decreasing, Abel summation gives
$\sum_{d=1}^{a}(N-d)(\varepsilon_d-\eta)\ge0$, whence
\[
\bar\sigma=\frac2N\sum_{d=1}^{a}(N-d)\psi_k(d)
=\frac{a(a+2)}{3k}+\frac2N\sum_{d=1}^{a}(N-d)\varepsilon_d
\ \ge\ \frac{a(a+2)}{3k}+a\eta .
\]
Combined with $\lambda>\rho(P)\ge\bar\sigma$ and \eqref{eq:gammalong} this gives
$\lambda-\gamma>g$.

For \eqref{eq:longnormgap} put $N=a+1$, so $a(a+2)=N^2-1$ and $mk+1\le N\le(m+1)k$.
Multiplying the difference of the two sides of \eqref{eq:longnormgap} by $3k^2$ gives
\[
F(N)=(2mk-N)(N^2-1)-3k^2N+3mk(k-1).
\]
On the interval $[mk+1,(m+1)k]$ we have $F''(N)=4mk-6N\le-2mk-6<0$, so $F$ is concave and it
suffices to check the endpoints. At $N=(m+1)k$,
\[
\frac{F((m+1)k)}k=k^2(m+1)(m^2-4)+m(3k-4)+1>0,
\]
since for $m=2$ this equals $6k-7>0$ and for $m\ge3$ all terms are nonnegative with the last
two positive. At $N=mk+1$, using $k\ge2$ and $m^2-3\ge1$,
\[
\begin{aligned}
\frac{F(mk+1)}k
&=k^2m(m^2-3)+k(m^2+3m-3)-5m\\
&\ge 4m+2(m^2+3m-3)-5m
 =2m^2+5m-6\\
&\ge 12>0 .
\end{aligned}
\]
Hence $F(N)>0$ throughout, which is \eqref{eq:longnormgap}.
\end{proof}

\begin{corollary}\label{cor:resolvent}
Let $\lambda>\rho(P)$ and $K=(\lambda I-P)^{-1}$. If $k\ge2$ and $a\ge k$, then for all
integers $j,h\ge1$,
\begin{equation}\label{eq:key}
\left|\tfrac12 q_j^{\tr}Kq_h\right|\ <\ \flo{\frac ak}\ \le\ \psi_k(a+j+h)-\psi_k(j+h).
\end{equation}
If $k=1$ and $a\ge1$, then $\tfrac12q_j^{\tr}Kq_h<\frac{a+1}2\le a=\psi_1(a+j+h)-\psi_1(j+h)$.
\end{corollary}

\begin{proof}
Let $k\le a<2k$, so $\flo{a/k}=1$. By \eqref{eq:resolvent}, Lemma~\ref{lem:middle} and
$\lambda-\gamma>g$,
\[
\tfrac12 q_j^{\tr}Kq_h\le\frac{\|q_j\|\|q_h\|}{2(\lambda-\gamma)}
<\frac{\max\{N,4(r+1)\}}{2g}\le1 .
\]
Let $a\ge2k$. By \eqref{eq:resolvent}, Lemmas~\ref{lem:norm} and \ref{lem:long},
\[
\tfrac12 q_j^{\tr}Kq_h\le\frac{\|q_j\|\|q_h\|}{2(\lambda-\gamma)}
<\frac{(a+1)\bigl(1+\frac{a(a+2)}{3k^2}\bigr)}{2g}<m .
\]
The second inequality in \eqref{eq:key} is Lemma~\ref{lem:ceiling}(2).

If $k=1$ then $q_t(i)=a-2i$ for every $t$, so
$\|q_t\|^2=\frac{(a+1)a(a+2)}3$ and $\bar\sigma=\frac{a(a+2)}3$; also $P=D$, so
\eqref{eq:pathnegative} gives $\gamma\le0$. Hence
\[
\tfrac12q_j^{\tr}Kq_h\le\frac{\|q_j\|\|q_h\|}{2\lambda}
<\frac{(a+1)a(a+2)/3}{2\bar\sigma}=\frac{a+1}2\le a\qquad(a\ge1).\qedhere
\]
\end{proof}

\subsection{The whole-pendant-path grafting lemma}

\begin{lemma}\label{lem:graft}
Let $T^\star$ be obtained from $T$ by the move \rm(M2) of Definition~\ref{def:moves}, with
$a\ge b\ge1$. Then for every integer $k\ge1$,
\begin{equation}\label{eq:graft}
\muk(T^\star)\ \ge\ \muk(T).
\end{equation}
If $C=\{v\}$, equality holds, both trees being paths. If $C\ne\{v\}$, then
\begin{equation}\label{eq:graftstrict}
\muk(T^\star)>\muk(T)\quad\Longleftrightarrow\quad H+a+b>k ,
\end{equation}
and when $H+a+b\le k$ the matrices $D_k(T^\star)$ and $D_k(T)$ coincide in the original
vertex labels.
\end{lemma}

\begin{proof}
If $C=\{v\}$ then $T$ and $T^\star$ are paths of order $a+b+1$ and \eqref{eq:graft} is an
equality. So assume $U\ne\emptyset$; write $\psi=\psi_k$, $x_0=v$,
$X=\{x_0,\dots,x_a\}$, $Y=\{y_1,\dots,y_b\}$ and $h_u=\dist_C(u,v)\ge1$ for $u\in U$.

\paragraph{Case 1: $a+b\le k$.}
Any two distinct vertices of $X\cup Y$ are at distance at most $a+b\le k$ in $T$, and in
$T^\star$ the set $X\cup Y$ induces a path of length $a+b$, so again all these distances are
at most $k$; thus all corresponding entries of $D_k$ equal $1$ in both trees. Distances
inside $U$, and from $U$ to $X$, are unchanged; distances from $u\in U$ to $y_j$ change from
$h_u+j$ to $h_u+a+j$. Hence $D_k(T^\star)\ge D_k(T)$ entrywise, proving \eqref{eq:graft}.
If $H+a+b\le k$, all mixed entries stay equal to $1$ and the matrices coincide. If
$H+a+b>k$, put $h=k-a-b+1$; then $1\le h\le H$, and since $C$ is a connected rooted tree of
height $H$ it contains a vertex $u$ of depth exactly $h$. Its distance to $y_b$ changes from
$h+b=k-a+1\le k$ to $h+a+b=k+1$, so the corresponding entry increases from $1$ to $2$.
Lemma~\ref{lem:PF}(2) gives strict inequality in \eqref{eq:graft}. This proves
\eqref{eq:graftstrict} in Case 1.

From now on $a+b>k$ (which forces $H+a+b>k$); it remains to prove strictness.

\paragraph{Block reduction.} In the order $X,Y,U$ set
\[
P_{i\ell}=\psi(|i-\ell|),\qquad \Sigma_{ij}=\psi(i+j),\qquad E_{iu}=\psi(i+h_u),
\qquad W=[\,\Sigma\ \ E\,],
\]
and let $F$ be the distance block on $Y\cup U$, i.e.\ $F_{y_jy_l}=\psi(|j-l|)$,
$F_{y_ju}=\psi(j+h_u)$, $F_{uu'}=\psi(\dist_C(u,u'))$. Then
\[
M=D_k(T)=\begin{pmatrix}P&W\\ W^{\tr}&F\end{pmatrix},
\qquad
M^\star=D_k(T^\star)=\begin{pmatrix}P&W'\\ W'^{\tr}&F'\end{pmatrix},
\]
where $W'=[\,R\Sigma\ \ E\,]$, because $\dist_{T^\star}(x_i,y_j)=a-i+j$, and $F'$ differs
from $F$ only in the $Y$--$U$ entries, which become $\psi(a+j+h_u)$.

Let $\lambda=\rho(M)$ and let $f=(\xi,\zeta)>0$ be a Perron vector, $\xi$ indexed by $X$ and
$\zeta$ by $Y\cup U$. By Lemma~\ref{lem:principal}, $\lambda>\rho(P)$, so
$K=(\lambda I-P)^{-1}$ exists and
\begin{equation}\label{eq:schur}
\xi=KW\zeta,\qquad (F+W^{\tr}KW)\zeta=\lambda\zeta .
\end{equation}
Since $P$, and hence $K$, commutes with $R$ and $R^{\tr}R=I$, we get
$(R\Sigma)^{\tr}K(R\Sigma)=\Sigma^{\tr}K\Sigma$: the $Y$--$Y$ and $U$--$U$ diagonal blocks of
$F+W^{\tr}KW$ and $F'+W'^{\tr}KW'$ agree. Their mixed $(y_j,u)$ entries differ by
\begin{equation}\label{eq:delta}
\delta_{j,h_u}
=\underbrace{\psi(a+j+h_u)-\psi(j+h_u)}_{=:\ \Delta_{j,h_u}}
-c_j^{\tr}(I-R)Kc_{h_u}
=\Delta_{j,h_u}-\tfrac12q_j^{\tr}Kq_{h_u},
\end{equation}
where we used $(R-I)K(R-I)=2K-RK-KR=2(I-R)K$ and $q_t=(R-I)c_t$.

For the trial vector $\hat f=(KW'\zeta,\zeta)\ne0$, substituting
$(P-\lambda I)KW'\zeta=-W'\zeta$ gives the exact identity
\begin{equation}\label{eq:rayleigh2}
\hat f^{\tr}(M^\star-\lambda I)\hat f
=\zeta^{\tr}\bigl(F'+W'^{\tr}KW'-\lambda I\bigr)\zeta
=2\sum_{u\in U}f_u\sum_{j=1}^{b}f_{y_j}\,\delta_{j,h_u},
\end{equation}
the last equality by \eqref{eq:schur} and the block structure just described. By
\eqref{eq:rayleigh}, positivity of the right-hand side implies $\muk(T^\star)>\lambda$.

\paragraph{Case 2: $a\ge k$, or $k=1$.}
Corollary~\ref{cor:resolvent} gives $\tfrac12q_j^{\tr}Kq_{h_u}<\Delta_{j,h_u}$, i.e.\
$\delta_{j,h_u}>0$ for all $j,u$. As $f>0$ and $Y,U\ne\emptyset$, the right side of
\eqref{eq:rayleigh2} is positive.

\paragraph{Case 3: $a<k<a+b$.}
Now all distances inside $X$ are less than $k$, so $P=J_N-I_N$ and $K$ acts on $\mathcal A$
as $(\lambda+1)^{-1}I$. Thus
\begin{equation}\label{eq:shortdelta}
\delta_{j,h}=\Delta_{j,h}-\frac{q_j^{\tr}q_h}{2(\lambda+1)} .
\end{equation}
Individual entries may be negative; we compare weighted column sums.

Let $S=\one^{\tr}f$. Since $a,b<k$, all $\psi$-values on $X\cup Y$ are $1$ or $2$. The Perron
equation at $v$ and at $y_j$ read
\[
\lambda f_v=\sum_{i=1}^{a}f_{x_i}+\sum_{l=1}^{b}f_{y_l}+\sum_{u}\psi(h_u)f_u,
\qquad
\lambda f_{y_j}=\sum_{i=0}^{a}\psi(i+j)f_{x_i}+\sum_{l\ne j}f_{y_l}
+\sum_u\psi(j+h_u)f_u .
\]
The first gives
\begin{equation}\label{eq:root}
(\lambda+1)f_v=S+\sum_{u\in U}\bigl(\psi(h_u)-1\bigr)f_u\ \ge\ S .
\end{equation}
Subtracting the first from the second and adding $f_{y_j}-f_v$ yields
\[
(\lambda+1)(f_{y_j}-f_v)
=\sum_{i=0}^{a}\bigl(\psi(i+j)-1\bigr)f_{x_i}
+\sum_{u\in U}\bigl(\psi(h_u+j)-\psi(h_u)\bigr)f_u .
\]
Every coefficient here is $0$ or $1$: indeed $i+j\le a+b\le2a<2k$ and $j\le b<k$
(Lemma~\ref{lem:ceiling}(2)). Hence the right-hand side lies in $[0,S]$, and with
\eqref{eq:root},
\begin{equation}\label{eq:perronbounds}
f_v\ \le\ f_{y_j}\ \le\ 2f_v\qquad(1\le j\le b).
\end{equation}

Put $d=k-a\ge1$ and $s=b-d=a+b-k>0$. For $j\le d$ we have $1\le a-i+j\le k$ and
$1\le i+j\le k$ for all $0\le i\le a$, so $q_j=0$; only the $s$ indices $j\in\{d+1,\dots,b\}$
can contribute. Also $|q_t(i)|\le1$ for all $t$ (base case of Lemma~\ref{lem:norm}), so
$q_j^{\tr}q_h\le\|q_j\|\|q_h\|\le N$. Finally $\Delta_{j,h}\in\{0,1\}$, and by
Lemma~\ref{lem:ceiling}(2) it vanishes exactly when $(j+h)\bmod k\in\{1,\dots,d\}$; the $b<k$
consecutive values $j=1,\dots,b$ have distinct residues, so at most $d$ increments vanish and
at least $s$ equal $1$. Combining with \eqref{eq:perronbounds},
\[
\sum_{j=1}^{b}f_{y_j}\delta_{j,h}
\ \ge\ s f_v-\frac{2f_v\,sN}{2(\lambda+1)}
= s f_v\Bigl(1-\frac N{\lambda+1}\Bigr)\ >\ 0 ,
\]
because $\lambda>\rho(P)=N-1$. Thus every weighted column sum in \eqref{eq:rayleigh2} is
positive, and $\muk(T^\star)>\lambda$. Cases 1--3 are exhaustive.
\end{proof}

\begin{remark}\label{rem:leafcounter}
The move (M1) is \emph{not} monotone for the distance spectral radius, so
Lemma~\ref{lem:hanging} has no distance analogue and Lemma~\ref{lem:graft} really needs the
whole-path graft. Let $T=S(2,2,1)$ and $T'=S(3,1,1)$, obtained from $T$ by moving the
terminal leaf of one $2$-leg to the end of the other $2$-leg; both have order $6$. In the
vertex order $(v,x_1,x_2,y_1,y_2,w)$ (centre, long leg, short leg, extra leaf) with $k=2$,
\[
D_2(T)=\begin{pmatrix}
0&1&1&1&1&1\\ 1&0&1&1&2&1\\ 1&1&0&2&2&2\\ 1&1&2&0&1&1\\ 1&2&2&1&0&2\\ 1&1&2&1&2&0
\end{pmatrix},
\quad
D_2(T')=\begin{pmatrix}
0&1&1&1&2&1\\ 1&0&1&1&1&1\\ 1&1&0&2&1&2\\ 1&1&2&0&2&1\\ 2&1&1&2&0&2\\ 1&1&2&1&2&0
\end{pmatrix},
\]
with characteristic polynomials
\[
(x+1)(x^2+3x+1)(x^3-4x^2-18x-9)
\quad\text{and}\quad
(x+1)\bigl(x^5-x^4-29x^3-63x^2-42x-8\bigr),
\]
respectively. The cubic $c(x)=x^3-4x^2-18x-9$ satisfies $c(34/5)=-241/125<0$ and
$c(x)\to\infty$, so its largest root exceeds $34/5$; the remaining factors of the first
polynomial have only negative roots, whence $\mu_2(T)>34/5$. The quintic
$u(x)=x^5-x^4-29x^3-63x^2-42x-8$ satisfies $u(6)=-2312<0$, $u(34/5)=237344/3125>0$, and is
strictly increasing on $[6,\infty)$ (indeed $u'(6)=1686>0$, $u''(6)=2718>0$ and
$u'''(x)=60x^2-24x-174>0$ for $x\ge6$); hence its largest root lies in $(6,34/5)$ and
$\mu_2(T')<34/5$. Therefore
\[
\mu_2(T')\ <\ \tfrac{34}5\ <\ \mu_2(T),
\qquad
\mu_2(T')\approx6.7830235,\quad \mu_2(T)\approx6.8288648 .
\]
\end{remark}

\subsection{Trees, and the proof of Theorem~\ref{thm:dist}}

\begin{proposition}\label{prop:trees-dist}
Let $T$ be a tree of order $n$ and $k\ge1$. Then $\muk(T)\le\muk(P_n)$. If $n\ge k+2$,
equality holds if and only if $T\cong P_n$.
\end{proposition}

\begin{proof}
If $T$ is not a path, root it at a leaf and choose a branching vertex $v$ farthest from the
root; every component of $T-v$ away from the root is a pendant path, and there are at least
two. Apply Lemma~\ref{lem:graft} to two of them, labelled so that $a\ge b$. The graft
decreases the number of leaves by exactly one: $x_a$ acquires degree $2$, $y_b$ stays a leaf,
$v$ keeps degree at least $2$, and no new leaf appears. By \eqref{eq:graft} the value of
$\muk$ does not decrease. Iterating, we reach a tree with two leaves, i.e.\ $P_n$, which
proves the weak inequality.

Assume $n\ge k+2$ and $T\not\cong P_n$. Since each graft removes exactly one leaf, the
sequence passes through a tree with exactly three leaves, i.e.\ a spider $S(a,b,c)$ with
$a,b,c\ge1$ and $a+b+c=n-1$, immediately before the last graft. For that graft the rooted
component $C$ is the third leg, so $H=c$ and
\[
H+a+b=n-1\ >\ k .
\]
By \eqref{eq:graftstrict} the last graft is strict, whence $\muk(T)<\muk(P_n)$.
\end{proof}

\begin{proof}[Proof of Theorem~\ref{thm:dist}]
If $n\le k+1$ then $G^k=K_n=P_n^k$ by Corollary~\ref{cor:powerfacts}(2), so all connected
graphs attain equality. Assume $n\ge k+2$ and let $T$ be a spanning tree of $G$. By
Corollary~\ref{cor:powerfacts}(3) and Proposition~\ref{prop:trees-dist},
\[
\muk(G)\ \le\ \muk(T)\ \le\ \muk(P_n).
\]
Suppose $G\not\cong P_n$. If $G$ has a vertex of degree at least three, extend three edges
incident with it to a spanning tree $T$ (any forest in a connected graph extends to a
spanning tree); then $T$ is not a path, so $\muk(T)<\muk(P_n)$ and we are done. Otherwise $G$
is a connected non-path graph of maximum degree at most two, i.e.\ a cycle. Deleting one
edge leaves a spanning path $T$; the ends of the deleted edge are at distance $1$ in $G$ and
$n-1>k$ in $T$, so $D_k(G)\le D_k(T)$ with strict inequality in at least one entry, and
Lemma~\ref{lem:PF}(2) gives $\muk(G)<\muk(T)\le\muk(P_n)$. In both cases equality fails;
conversely $G\cong P_n$ gives equality.
\end{proof}

\section{Concluding remarks and open problems}\label{sec:conc}

The two main theorems identify the same extremal graph for two different spectral quantities:
for each fixed $k$, the power of the path minimizes the adjacency spectral radius and
maximizes the distance spectral radius among the $k$-th powers of connected graphs of order
$n$. The equality cases change only at the natural threshold $n=k+2$, below which every
$k$-th power is complete.

The proofs use different local transformations. For adjacency spectral radius, transferring
a single terminal vertex between two pendant paths is sufficient. For distance spectral
radius, this move can have the wrong monotonicity, and the entire shorter pendant path must
be transferred. The comparison for the latter move is based on a Schur complement and on
bounds for the largest eigenvalue of the antisymmetric restriction of a path-power distance
matrix. The examples below illustrate why the adjacency rearrangement argument cannot be
applied unchanged to distance matrices.

\begin{table}[h]
\centering
\begin{tabular}{lll}
\toprule
Operation & adjacency spectral radius & distance spectral radius\\
\midrule
(M1) terminal-vertex transfer & non-increasing (Lemma~\ref{lem:hanging}) & not monotone\\
(M2) whole-path graft & non-increasing (Corollary~\ref{cor:graft-adj}) & non-decreasing (Lemma~\ref{lem:graft})\\
\bottomrule
\end{tabular}
\caption{Effect of the two transformations on the spectral radii of graph powers.}
\label{tab:ops}
\end{table}

The terminal-vertex transfer is not monotone for the distance spectral radius. Let
$T=S(2,2,1)$ and let $T'=S(3,1,1)$ be obtained by transferring the terminal vertex of one
leg of length $2$ to the end of the other leg of length $2$. Both trees have order $6$. In
the vertex order $(v,x_1,x_2,y_1,y_2,w)$, with $v$ the centre, and for $k=2$, their distance
matrices are
\[
D_2(T)=\begin{pmatrix}
0&1&1&1&1&1\\ 1&0&1&1&2&1\\ 1&1&0&2&2&2\\ 1&1&2&0&1&1\\ 1&2&2&1&0&2\\ 1&1&2&1&2&0
\end{pmatrix},\qquad
D_2(T')=\begin{pmatrix}
0&1&1&1&2&1\\ 1&0&1&1&1&1\\ 1&1&0&2&1&2\\ 1&1&2&0&2&1\\ 2&1&1&2&0&2\\ 1&1&2&1&2&0
\end{pmatrix}.
\]
Their characteristic polynomials are, respectively,
\[
(x+1)(x^2+3x+1)(x^3-4x^2-18x-9)
\quad\text{and}\quad
(x+1)(x^5-x^4-29x^3-63x^2-42x-8).
\]
The cubic $c(x)=x^3-4x^2-18x-9$ satisfies $c(34/5)=-241/125<0$ and tends to $+\infty$;
its largest root is therefore greater than $34/5$. The other two roots of the first
characteristic polynomial are negative. For the quintic
$u(x)=x^5-x^4-29x^3-63x^2-42x-8$, we have $u(6)=-2312$ and
$u(34/5)=237344/3125$. Moreover,
$u'(6)=1686$, $u''(6)=2718$, and
$u^{\prime\prime\prime}(x)=60x^2-24x-174>0$ for $x\geq6$. Hence $u$ is strictly increasing on $[6,\infty)$,
and its largest root lies in $(6,34/5)$. It follows that
\[
\mu_2(T')<\frac{34}{5}<\mu_2(T),
\qquad
\mu_2(T')\approx6.7830235,\quad \mu_2(T)\approx6.8288648.
\]

The reflection compression used for adjacency matrices also fails for distance matrices. Let
$K=D(P_5^2)$ and take
\[
f=(2,1,1,2,1)^{\tr},\qquad z=(1,1,1,2,2)^{\tr}.
\]
The vector $z$ is obtained from $f$ by placing the larger entry of each reflection pair on
the right. Although $\|f\|^2=\|z\|^2=11$, direct calculation gives
\[
f^{\tr}Kf=52>50=z^{\tr}Kz.
\]
Thus monotonicity of the distance kernel alone does not provide the rearrangement inequality
needed in the adjacency proof.

Several questions remain open.

\begin{question}
Determine the graphs minimizing or maximizing the largest eigenvalue of
\[
A_\alpha(G^k)=\alpha\operatorname{Deg}(G^k)+(1-\alpha)A(G^k),
\qquad 0\leq\alpha<1,
\]
among connected graphs of order $n$. For $k=2$, the extremal problem is settled by
\cite{ChenLiWangMeng2019}. For $k\geq3$, is the path extremal for small $\alpha$, and how
does the answer depend on $\alpha$?
\end{question}

\begin{question}
Find analogues of the two main theorems for the signless Laplacian and Laplacian spectral
radii of $G^k$. In particular, does the whole-path graft (M2) preserve a useful monotonicity
property for either matrix?
\end{question}

\begin{question}
Determine the extremal trees for $\rho_k$ and $\mu_k$ within classes with fixed maximum
degree, diameter, or matching number. Such constraints are of particular interest for the
distance spectral radius, since the $k$-th power of a star is complete whenever $k\geq2$.
\end{question}

\begin{question}
Beyond the largest eigenvalue, determine extremal graphs for the second largest adjacency
eigenvalue or the least distance eigenvalue of $G^k$, over connected graphs or trees of order
$n$.
\end{question}

\begin{question}
For which nonnegative nondecreasing functions $f:\mathbb Z_{\ge0}\to[0,\infty)$ with
$f(0)=0$ is the path extremal for the largest eigenvalue of the matrix
$[f(\dist_G(u,v))]_{u,v\in V(G)}$ among connected graphs of order $n$? The two main theorems
correspond to the kernels $f_k(t)=\mathbf 1[1\le t\le k]$ and
$f(t)=\lceil t/k\rceil$, respectively. The counterexample above shows that monotonicity
of $f$ alone does not guarantee the rearrangement inequality used in the adjacency argument.
\end{question}

\section*{Declaration on the use of AI}
The authors used generative AI tools to assist in discussing proof strategies, checking proofs, and improving the exposition. The authors take full responsibility for the mathematical arguments, results, and conclusions, all of which were carefully reviewed and verified by them.

\section*{Declaration of competing interest}

The authors declare that they have no known competing financial interests or personal relationships that could have appeared to influence the work reported in this article.

\section*{Declaration of Funding}

The authors declare that they have not received any external funding for the work done in this manuscript.

\section*{Acknowledgments}
The first author Joyentanuj Das acknowledge Hitesh Kumar, for introducing the problem and helping with proof ideas.

\section*{Data availability}

Data sharing is not applicable to this article as no datasets were generated or analyzed during the current study.

\end{document}